\documentclass[11pt]{amsart}
\usepackage{amsmath,amssymb}
\usepackage{setspace,xypic}
\usepackage[all]{xy}

\newcommand\inv{^{-1}}

\newcommand\fh{\mathcal H}

\newcommand\G{\Gamma}

\newcommand\g{\gamma}
\newcommand\Ta{\mathbb T}

\newcommand\Za{\mathbb Z}

 \DeclareMathOperator{\Isom}{Isom}

 \DeclareMathOperator{\Hom}{Hom}
\DeclareMathOperator{\Aut}{Aut}

\DeclareMathOperator{\Out}{Out}
\DeclareMathOperator{\Trace}{Trace}\DeclareMathOperator{\gap}{gap}
\DeclareMathOperator{\pgap}{pgap}

\newtheorem{theorem}{Theorem}[section]

\newtheorem{lemma}[theorem]{Lemma}

\newtheorem{conjecture}[theorem]{Conjecture}

\begin{document}

\title[$\Out(F_n)$ and the spectral gap]{$\Out(F_n)$ and the spectral gap
conjecture}
\author{David Fisher}
\thanks{Author partially supported by NSF grant DMS-0226121 and DMS-0541917.}

\begin{abstract}
For $n>2$, given $\phi_1,\ldots,\phi_n$ randomly chosen isometries
of $S^2$, it is well-known that the group $\G$ generated by
$\phi_1,\ldots,\phi_n$ acts ergodically on $S^2$.  It is
conjectured in \cite{GJS} that for almost every choice of
$\phi_1,\ldots,\phi_n$ this action is {\em strongly ergodic}. This
is equivalent to the spectrum of
$\phi_1+\phi_1{\inv}+{\cdots}+\phi_n+\phi_n^{\inv}$ as an operator
on $L^2(S^2)$ having a spectral gap, i.e. all eigenvalues but the
largest one being bounded above by some $\lambda_1<2n$.  (The
largest eigenvalue $\lambda_0$, corresponding to constant
functions, is $2n$.)

In this article we show that if $n>2$, then either the conjecture
is true or almost every $n$-tuple fails to have a gap. In fact,
the same result is holds for any $n$-tuple $\phi_1,\ldots, \phi_n$
in any any compact group $K$ that is an almost direct product of
$SU(2)$ factors with $L^2(S^2)$ replaced by $L^2(X)$ where $X$ is
any homogeneous $K$ space. A weaker result is proven for $n=2$ and
some conditional results for similar actions of $F_n$ on
homogeneous spaces for more general compact groups.
\end{abstract}

\maketitle

\section{Introduction.}
\label{section:results}

Let $\phi_1,\ldots,\phi_n$ be any finite collection of elements of
$SU(2)$ and let $L_0^2(SU(2))$ be the orthogonal complement of the
constant functions in $L^2(SU(2))$.  The operator
$\phi_1+\phi_1{\inv}+{\cdots}+\phi_n+\phi_n^{\inv}$ is a
self-adjoint operator on $L_0^2(SU(2))$ and has discrete spectrum
which is a subset of $\mathbb R^+$.  Let $\lambda_1$ be the
supremum of the eigenvalues for this operator. It is clear that
$\lambda_1{\leq}2n$. The following is conjectured in \cite{GJS}:

\begin{conjecture}
\label{conjecture:spectralgap} For $n{\geq}2$ and almost every
collection $\phi_1,\ldots,\phi_n$, we have $\lambda_1<2n$ for
$\phi_1+\phi_1{\inv}+{\cdots}+\phi_n+\phi_n^{\inv}$.
\end{conjecture}

\noindent This conjecture is referred to as the spectral gap
conjecture and is a question in \cite{LPS}.   The conjecture is
only known for $n$-tuples which have, up to conjugacy in $SU(2)$,
all matrix entries of all $\phi_i$ algebraic. This is a recent
result of Bourgain and Gamburd building on earlier work of
Gamburd, Jakobson and Sarnak \cite{BG,GJS}. This set of $n$-tuples
for which the conjecture is known has zero measure. See
\cite[Theorem 3.2]{F} and \cite{KR} and the references there for
weaker related results.

The main result of this paper is the following:

\begin{theorem}
\label{theorem:01law} Assume $n{\geq}3$.  Then either
$\lambda_1<2n$ for almost every $\phi_1,\ldots,\phi_n$   or
$\lambda_1=2n$ for almost every $\phi_1,\ldots,\phi_n$.
\end{theorem}

\noindent In particular, by Theorem \ref{theorem:01law}, to prove
Conjecture \ref{conjecture:spectralgap} it suffices to establish a
spectral gap for any set of positive measure in $SU(2)^n$.

It is well-known that for almost every $\phi_1,\ldots,\phi_n$, the
group generated by $\phi_1,\ldots,\phi_n$ is a free group on $n$
generators, $F_n$. The space of $n$-tuples $\phi_1,\ldots,\phi_n$
can be parametrized as $\Hom(F_n,SU(2))$. The main new ingredient
in Theorem \ref{theorem:01law} is the use of symmetries of $F_n$
and in particular ergodicity of a the action of $\Aut(F_n)$ on
$\Hom(F_n,SU(2))$ where $\Aut(F_n)$ is the automorphism group of
$F_n$.  

Theorem \ref{theorem:01law} remains true when $SU(2)$ is replaced
by any compact Lie group $K$ which is an almost direct product of
copies of $SU(2)$ and $SU(1)$. For any compact Lie group with
$SU(1)=S^1$ factors, analogues of Theorem \ref{theorem:01law} are
not interesting, as it is easy to see in that case that there is
no spectral gap on a set of full measure.

The interest in spectral gaps for finite collections of elements
in $SU(2)$ originally derives from the Banach-Ruscewiecz
conjecture and its proof. This states that, for $m>1$, the unique
finitely additive rotationally invariant measure on $S^m$ is the
Haar measure. Rosenblatt showed that this was equivalent to
finding a finite subset $\phi_1,{\ldots},\phi_n$ in $\Isom(S^m)$
with a spectral gap for the action on $L^2_0(S^m)$ \cite{Ro}. For
$n>1$, a spectral gap for $\phi_1,{\ldots}\phi_n$ in $\Isom(S^m)$
on $L^2_0(S^m)$ is easily seen to be equivalent to a spectral gap
for $\phi_1,{\ldots}\phi_n$ on $L^2_0(\Isom(S^m))$ since the same
representations of $\Isom(S^m)$ occur in $L^2_0(S^m)$ and
$L^2_0(\Isom(S^m))$, just with different multiplicities.
Similarly, Conjecture \ref{conjecture:spectralgap} is equivalent
to the same conjecture with $SO(3)$ in place of $SU(2)$.  For
$n>3$ Sullivan and Margulis independently exhibited such subsets
with a spectral gap, each by finding a homomorphism from a group
$\G$ with property $(T)$ of Kazhdan to $\Isom(S^m)$ \cite{Ma,Su}.
For $n=2,3$, subsets of $\Isom(S^m)$ with a spectral gap were
first exhibited by Drinfeld using methods of automorphic forms
\cite{Dr}.  Later work on the subject was motivated by the fact
that if $\phi_1,{\ldots},\phi_n$ have a spectral gap, then the
orbits under the resulting action of $F_n$ on $S^m$ equidistribute
with exponential speed.  In \cite{LPS}, the authors show how to
find $\phi_1,{\ldots}\phi_n$ with optimal equidistribution
properties, again using deep results on automorphic forms.  In
\cite{GJS}, the authors prove the existence of
$\phi_1,{\ldots}\phi_n$ in $SU(2)$ with a spectral gap without
using heavy machinery from the theory of automorphic forms and
also discuss several related issues.  For more discussion see
\cite{GJS,Lu,Sa}.

As mentioned above, the key step in the proof of all results here
is to use the ergodic theory of the action of $\Aut(F_n)$ on
$\Hom(F_n,K)$. In fact, since it is easy to check that the
spectral gap is invariant under conjugation in $SU(2)$ it is
easier to work with the action of $\Out(F_n)$ on
$\Hom(F_n,SU(2))/SU(2)$ instead. The group $\Aut(F_n)$ of
automorphisms of $F_n$ acts on $\Hom(F_n,K)$ and this action
descends to an action of the outer automorphism group $\Out(F_n)$
on $\Hom(F_n,K)/K$. The $\Aut(F_n)$ action preserves the measure
on $Hom(F_n,K)$ given by identifying this space with $K^n$ and
taking Haar measure.  The $\Out(F_n)$ action preserves the measure
on $\Hom(F_n,K)/K$ given by realizing $\Hom(F_n,K)$ as $K^n$,
taking Haar measure on each factor, and dividing by the
conjugation action of $K$ to obtain the quotient $\Hom(F_n,K)/K$.
The dynamics of this action have received relatively little
attention, but the analogous action of the mapping class group on
$\Hom(S,K)/K$ where $S$ is the fundamental group of a surface, has
been studied more extensively, see the recent survey \cite{Go3}
which is also a good introduction to dynamics of group actions on
representation varieties. Essentially the only known result for
the action of $\Out(F_n)$ on $\Hom(F_n,K)/K$ is due to Goldman who
shows that the action is weakly mixing when $k{\geq}3$ and $K$ is
an almost direct product of $SU(2)$ and $SU(1)$ factors. This is
proven in \cite{Go2} using the main results of \cite{Go1}. The
remaining ingredient in the proof of Theorem \ref{theorem:01law}
is to construct a measurable function $f$ on
$\Hom(F_n,SU(2))/SU(2)$ that is $\Out(F_n)$ invariant and takes
the value $1$ for actions with a spectral gap and the value zero
for actions without a spectral gap. We will construct the function
$f$ in \S \ref{section:proofs}. In the next section, we state some
other variants of Theorem \ref{theorem:01law}. and recall some
results about the ergodic theory of actions on moduli spaces from
\cite{Go1,Go2,PX1,PX2}. In section  \S \ref{section:proofs} we
prove all of our results.

\section{Further results and group actions on representation
varieties.} \label{section:furtherresults}

A key ingredient in our proof of Theorem \ref{theorem:01law} is
the following result of Goldman.

\begin{theorem}[Goldman]
\label{theorem:bill} Let $K$ be a compact group which is an almost
direct product of $SU(2)$ and $SU(1)$ factors.  If $n>2$, then the
action of $\Out(F_n)$ on $\Hom(F_n,K)/K$ is ergodic.
\end{theorem}

When $n=2$, there are non-constant function on
$\Hom(F_2,SU(2))/SU(2)$ or even $\Hom(F_n,K)/K$ which are easily
seen to be $\Out(F_2)=SL(2,\Za)$ invariant.  For $K=SU(2)$, one
such function is simply $g(\rho)=\Trace([\rho(a),\rho(b)])$ where
$a,b$ are a basis for $F_2$. In this case invariance follows from
the fact that the set of commutators is $\Aut(F_2)$ invariant and
that trace is conjugation invariant. For general $K$, we need a
few facts before we can define an analogous function.  It is
well-known that every element of $K$ is contained in a maximal
torus $T<K$ and that all such tori are conjugate in $K$.  This
allows us to parametrize the conjugacy classes in $K$ as $T/W$
where $W<K$ is the Weyl group, i.e. the normalizer of $T$ divided
by the centralizer of $T$. Given $K$, we define
$g:\Hom(F_2,K)/K{\rightarrow}T/W$ by taking the representative of
the conjugacy class of $[\rho(a),\rho(b)]$. Again this function is
invariant, since the commutator is $\Out(F_2)$ invariant. The
following result of Pickrell and Xia is essentially \cite[Theorem
2.1.4]{PX1}.  In the case where $K$ is as in Theorem
\ref{theorem:bill}, the result is contained in \cite{Go1}.

\begin{theorem}[Pickrell-Xia]
\label{theorem:px} For any compact Lie group $K$, the map
$g:\Hom(F_2,K)/K{\rightarrow}T/W$ defined above is an ergodic
decomposition for the action of $\Out(F_2)$ on $\Hom(F_2,K)/K$.
\end{theorem}

\noindent In both \cite{Go1} and \cite{PX1}, $\Out(F_2)$ is
considered as the mapping class group of a once punctured torus.

By viewing $g$ as an ergodic decomposition, we are writing the
measure on $\Hom(F_2,K)/K$ as an integral over $T/W$ of measures
on the level sets of $g$. In fact, level sets of $g$ are
generically smooth submanifolds and these measures are smooth
measures. It is clear that we can view $g$ as a function on
$\Hom(F_2,K)$ instead. In the following result, $\lambda_1$ is
again the supremum of eigenvalues for the operator
$\rho(a)+\rho(b)+\rho(a){\inv}+\rho(b){\inv}$ on $L^2_0(K)$.

\begin{theorem}
\label{theorem:n=2} Let $a,b$ be a basis for $F_2$ and let $K$ be
any compact group. Let $X$ be the subset of $\Hom(F_2,K)$ such
that $\lambda_1<4$ for
$\rho(a)+\rho(b)+\rho(a){\inv}+\rho(b){\inv}$. Then for almost
every $a$ in the image of $g$ (with the pushforward measure), the
set $X{\cap}g{\inv}(a)$ has either zero measure or full measure in
$g{\inv}(a)$.
\end{theorem}

\noindent This theorem is proven exactly as Theorem
\ref{theorem:01law}, using Theorem \ref{theorem:px} in place of
Theorem \ref{theorem:bill}. The main obstruction to a variant of
Theorem \ref{theorem:01law} for general $K$ is the lack of an
analogue of the main result of \cite{Go2} for general $K$. The
following conditional result also follows from the proof of
Theorem \ref{theorem:01law}.

\begin{theorem}
\label{theorem:conditional} Let $a_1,a_2,\ldots,a_n$ be a basis
for $F_n$. Let $X$ be the subset of $\Hom(F_n,K))$ such that
$\lambda_1<2n$ for
$\rho(a_1)+\rho(a_1{\inv})+\cdots+\rho(a_n)+\rho(a_n){\inv}$. Then
if $n>2$, the measure of $X$ is either $0$ or $1$ provided the
$\Out(F_n)$ action on $\Hom(F_n, K)/K$ is ergodic.
\end{theorem}

\section{Proofs of the main results.} \label{section:proofs}

The proof of all results here depend on another, equivalent,
definition of the spectral gap. The following discussion and the
first two lemmas of this section are standard, but we include them
for completeness. We can define the {\em spectral gap} for a
unitary representation $\rho$ of a finitely generated group $\G$
with generating set $S$ on a Hilbert space $\fh$ to be the largest
$\varepsilon$ such that for each $v{\in}\fh$ there is some $\g$ in
$S$ such that:
$$\|v-\rho(\g)v\|{\geq}\varepsilon\|v\|.$$
\noindent Note that the spectral gap depends on the generating set
$S$. This is because a choice of generating sets determines a
particular basis of neighborhoods of the trivial representation in
the Fell topology. That having a non-zero spectral gap in this
sense is equivalent to the definition of spectral gap given above
is more or less immediate from the definition of the Fell topology
on the unitary dual of group, but we give a proof below for
completeness.The following standard lemma shows that having a
non-zero spectral gap is independent of generating set.

\begin{lemma}
\label{lemma:gapandgenerators} Let $\G$ be a finitely generated
group and let $S_1$ and $S_2$ be two generating sets for $\G$.
Then $\G$ has  non-zero spectral gap for $S_1$ if and only if it
has a non-zero spectral gap for $S_2$.
\end{lemma}

\begin{proof}
Let $\varepsilon$ be the spectral gap for $(\G,S_1)$.  Let $n$ be
the smallest integer such that every element of $S_1$ can be
written as a word of length $n$ in the generators $S_2$.  Then we
claim that the spectral gap for $(\G,S_2)$ is at least
$\frac{\varepsilon}{n}$.  If not we have a vector $v$ in a
representation $\sigma$ of $\G$ on a Hilbert space $\fh$ such that
$$\|\sigma(\g)v-v\|<\frac{\varepsilon}{n}$$
\noindent for all $\g$ in $S_2$.  Writing any $\tilde \g{\in}S_1$
as $\g_1{\cdots}\g_i$ where $i\leq{n}$ and using a standard
telescoping sum argument, this implies that
$$\|\sigma(\tilde \g)v-v\|<{\varepsilon}$$
\noindent a contradiction.  Reversing the roles of $S_1$ and $S_2$
completes the proof.
\end{proof}

The following standard lemma implies that our two definitions of
spectral gap are equivalent.

\begin{lemma}
\label{lemma:gaps} Let $\G$ be a finitely generated group with
generators $g_1,g_2,\ldots,g_n$ and $\rho$ a unitary
representation of $\G$ on a Hilbert space $\fh$. Then $\rho$ has a
spectral gap if and only if the norm the operator
$\rho(g_1)+\cdots+\rho(g_n){\inv}$ is strictly less than $2n$.
\end{lemma}

\begin{proof} If we adjoin the identity to any generating set as
$\g_{n+1}$, the norm of the operator
$\rho(g_1)+\cdots+\rho(g_{n+1}){\inv}$  is simply $2$ plus the
norm of the operator $\rho(g_1)+\cdots+\rho(g_n){\inv}$. Combined
with Lemma \ref{lemma:gapandgenerators} this means it suffices to
prove the current lemma for generating sets that contain the
identity. This reduces to the following elementary fact about
Hilbert spaces: given $k$ unit vectors $v_1,\ldots,v_k$ not all of
which are equal, $\|\frac{1}{k}\sum v_k\|<1-f(v_1,\ldots,v_k)$
where $f$ is a positive function of the diameter of the set
$v_1,\ldots, v_k$ which goes to zero only when the diameter goes
to zero. (It is not too hard to write down $f$ explicitly.)
\end{proof}

We now prove a lemma that suffices to prove all the theorems in
the previous section.

\begin{lemma}
\label{lemma:measurablegap} Fix a finitely generated group $\G$
and a generating set $\g_1,\ldots,\g_n$.  Define the function
$\gap(\rho)$ to be the spectral gap of $(\rho(\G),S)$ acting on
$L_0^2(K)$ where $K$ is compact Lie group and $\rho$ is in
$\Hom(\G,K)$. Then $\gap$ is a well-defined measurable function on
$\Hom(\G,K)/K$.
\end{lemma}

\begin{proof}
Fix a Riemannian metric on $K$ invariant under both left and right
multiplication.  Let $\Delta$ be the associated Laplacian. Recall
that $L_0^2(K)$ decomposes as a Hilbertian direct sum
$$\oplus_{\lambda}V_{\lambda}$$
\noindent where $\lambda$ runs over non-zero eigenvalues of
$\Delta$ and each $V_{\lambda}$ is a bi-$K$ invariant finite
dimensional space of smooth functions on $K$.  Let
$S(V_{\lambda})$ be the unit sphere in $V_{\lambda}$.  The action
of $K$ on $S(V_{\lambda})$ is smooth, so for every representation
$\rho:\G{\rightarrow}K$, we have a smooth action $\rho_{\lambda}$
of $\G$ on $S(V_{\lambda})$ and $\rho_{\lambda}$ depend smoothly
on $\rho$. Therefore the function $\rho(\g_i)v-v$ is a smooth
function on $\Hom(\G,K){\times}S(V_{\lambda})$. The function
$\|\rho(\g_i)v-v\|$ is continuous on
$\Hom(\G,K){\times}S(V_{\lambda})$ and so the function
$\gap_{\lambda}(\rho,\g_i)=\min_{S(V_{\lambda})}\|\rho(\g_i)v-v\|$
is continuous on $\Hom(\G,K)$.  Therefore:

$$\widetilde{\gap}(\rho)=\inf_{\lambda}\max_{\g_1,\ldots,\g_n}\gap_{\lambda}(\rho,\g_i)$$

\noindent is a measurable function on $\Hom(\G,K)$. It is
immediate from the definition that $\widetilde{\gap}(\rho)$ is in
fact the spectral gap for $(\rho,S)$ and that it is invariant
under conjugation.
\end{proof}

\noindent We now proceed to prove the theorems stated in the
introduction.

\begin{proof}[Proof of Theorems \ref{theorem:01law}, \ref{theorem:n=2}
and \ref{theorem:conditional}] We define a function $\pgap(\rho)$
on the space $\Hom(\G,K)/K$ such that $\pgap(\rho)=1$ if
$\gap(\rho)>0$ and $\pgap(\rho)=0$ otherwise. By Lemma
\ref{lemma:gapandgenerators} the function $\pgap$ is $\Out(F_n)$
invariant, an automorphism of $F_n$ simply changes the generating
set for which we want a gap.  By Lemma \ref{lemma:measurablegap}
the function $\pgap$ is measurable.

To complete the proof of Theorem \ref{theorem:01law}, we note that
by Theorem \ref{theorem:bill}, the action of $\Out(F_n)$ on
$\Hom(\G,K)/K$ is ergodic as long as $K$ is locally a product of
$SU(2)$ and $SU(1)$ factors and $n>2$.  This implies that $\pgap$
is either almost everywhere one or almost everywhere zero.
Similarly, to complete the proof of Theorem \ref{theorem:n=2}, we
recall that by Theorem \ref{theorem:px}, the level sets of the
function $g$ defined in \S \ref{section:furtherresults} are
ergodic components for the action of $\Out(F_2)$ on
$\Hom(F_2,K)/K$.  This immediately implies the statement of the
theorem.  The proof of Theorem \ref{theorem:conditional} is the
same as the proof of Theorem \ref{theorem:01law}.
\end{proof}

\section{Speculation and questions} \label{section:speculation}

Theorem \ref{theorem:01law} allows one to prove Conjecture
\ref{conjecture:spectralgap} by proving the existence of a
spectral gap on any set of positive measure in $\Hom(F_n,SU(2))$.
It also leaves one with the impression that the large group of
symmetries of $F_n$ might be relevant to a proof of Conjecture
\ref{conjecture:spectralgap}.

The following conjecture seems natural in the context of this
work:

\begin{conjecture}
\label{conjecture:stronglyergodic} The representation of
$\Out(F_n)$ on $L^2_0(\Hom(F_n,K)/K)$ has a spectral gap for
$n>3$.
\end{conjecture}

\noindent One can reformulate this as saying that the trivial
representation is isolated in the representation of $\Out(F_n)$ on
$L^2(\Hom(F_n,K)/K)$, in which case the conjecture also makes
sense for $n=2$.  It may be possible to prove Conjecture
\ref{conjecture:stronglyergodic} for $n$ large enough, using the
fact that $\Out(F_n)$ is generated by torsion elements, see e.g.
\cite{BV,Zu}, and an argument like the one given by Schmidt in
\cite{Sch} for strong ergodicity of the $SL(2,\Za)$ action on
$\Ta^2$.   In the case when $K$ is abelian, the conjecture is true
and originally due to Rosenblatt \cite{Ro}.  When $K$ is abelian,
stronger statements in this direction, including strong ergodicity
of many subgroups, follow from work of Furman and Shalom
\cite{FS}.  It is tempting to hope for some duality that links
Conjecture \ref{conjecture:stronglyergodic} to Conjecture
\ref{conjecture:spectralgap}, but this hope seems naive.  Any
attempt to link the two conjectures must take account of the fact
that Conjecture \ref{conjecture:stronglyergodic} is true when $K$
is abelian, and the analogue of Conjecture
\ref{conjecture:spectralgap} fails in that setting.

It is also worthwhile to compare this paper to work where
relations are sought between spectral gaps in certain (or all)
representations of $\Out(F_n)$ and expansion properties of various
families of finite groups, see particulary \cite{LP,GP}.  In
particular, it seems likely that a strong version of Conjecture
\ref{conjecture:spectralgap} should imply Conjecture
\ref{conjecture:stronglyergodic} via an argument similar to the
one in \cite{GP}.

\bigskip
\noindent Department of Mathematics, Indiana University, Rawles Hall,
Bloomington, IN 47401.\\


\begin{thebibliography}{Palais}




\bibitem[BV]{BV}
Bridson, Martin R.; Vogtmann, Karen Homomorphisms from
automorphism groups of free groups. {\it Bull. London Math. Soc.}
35  (2003),  no. 6, 785--792.


\bibitem[BG]{BG}
Bourgain, J. and Gamburd, A. On the spectral gap for finitely
generated subgroups of $SU(2)$, preprint 2006.

\bibitem[Dr]{Dr}
Drinfeld, V. G.
Finitely-additive measures on $S\sp{2}$ and $S\sp{3}$, invariant with respect to rotations. (Russian)
{\it Funktsional. Anal. i Prilozhen.} 18 (1984), no. 3, 77.

\bibitem[F]{F}
Fisher, David. First cohomology and local rigidity of group
actions, preprint, arxiv math.DG/0505520.

\bibitem[FS]{FS}
Furman, Alex; Shalom, Yehuda. Sharp ergodic theorems for group
actions and strong ergodicity. {\it Ergodic Theory Dynam. Systems}
19 (1999),  no. 4, 1037--1061.

\bibitem[GJS]{GJS}
Gamburd, Alex; Jakobson, Dmitry; Sarnak, Peter. Spectra of elements
in the group ring of ${\rm SU}(2)$. {\it J. Eur. Math. Soc. (JEMS)}
1 (1999),  no. 1, 51--85.

\bibitem[GP]{GP}
Gamburd, Alex; Pak, Igor.  Expansion of product replacement
graphs, to appear {\em Combinatorica}.

\bibitem[Go1]{Go1} Goldman, William M. Ergodic theory on moduli spaces.
{\it Ann. of Math. (2)}  146  (1997),  no. 3, 475--507.

\bibitem[Go2]{Go2} Goldman, William M. An ergodic action of the
outer automorphism group of a free group, preprint, arxiv
math.DG/0506401 .

\bibitem[Go3]{Go3}
Goldman, William M. Mapping Class Group Dynamics on Surface Group
Representations, preprint, arxiv math.GT/0509114.


\bibitem[KR]{KR}
Kaloshin, V.; Rodnianski, I. Diophantine properties of elements of
${\rm SO}(3)$. {\it Geom. Funct. Anal.}  11  (2001),  no. 5,
953--970.

\bibitem[Lu]{Lu}
A. Lubotzky. {\it Discrete groups, expanding graphs and invariant
measures.} Birkhauser, Basel 1994

\bibitem[LPS]{LPS} A. Lubotzky, R. Phillips, P. Sarnak.
Hecke operators and distributing points on $S^2$. I. {\it Comm.
Pure Appl. Math.} 39(S), S149--S186 (1986) II. {\it Comm. Pure
Appl. Math.} 40(4), 401--420 (1987).

\bibitem[LP]{LP}
Lubotzky, Alexander; Pak, Igor. The product replacement algorithm
and Kazhdan's property (T). {\it J. Amer. Math. Soc.} 14 (2001),
no. 2, 347--363

\bibitem[Ma]{Ma}
G. Margulis. Some remarks on invariant means, Monatschefte fur
Mathematik 90 (1980), 233--235.


\bibitem[PX1]{PX1} Pickrell, Doug; Xia, Eugene Z. Ergodicity of mapping
class group actions on representation varieties. I. Closed
surfaces. {\it Comment. Math. Helv.}  77  (2002),  no. 2,
339--362.

\bibitem[PX2]{PX2}
Pickrell, Doug; Xia, Eugene Z.
Ergodicity of mapping class group actions on representation
varieties. II. Surfaces with boundary. (English. English summary)
{\it Transform. Groups} 8 (2003), no. 4, 397--402.


\bibitem[Ro]{Ro}
J. Rosenblatt. Uniqueness of invariant means for measure preserving
transformations, Trans. AMS 265 (1981), 623--636.

\bibitem[Sa]{Sa}
Sarnak, Peter. {\it Some applications of modular forms.} Cambridge
Tracts in Mathematics, 99. Cambridge University Press, Cambridge,
1990.

\bibitem[Sch]{Sch}
Schmidt, Klaus. Asymptotically invariant sequences and an action of
${\rm SL}(2,\,Z)$ on the $2$-sphere.  Israel J. Math.  37  (1980),
no. 3, 193--208.

\bibitem[Su]{Su}
D. Sullivan. For $n > 3$ there is only one finitely additive
rotationally invariant measure on the $n$-sphere on all Lebesgue
measurable sets, Bull. AMS 1 (1981), 121--123.

\bibitem[Zu]{Zu}
Zucca, Paola. On the $(2,2\times2)$-generation of the automorphism
groups of free groups.  {\it Istit. Lombardo Accad. Sci. Lett.
Rend.} A 131  (1997),  no. 1-2, 179--188 (1998).


\end{thebibliography}
\end{document}